\documentclass[11pt]{article}
\usepackage{latexsym, amsfonts}
\newcommand{\be}{\begin{equation}}
\newcommand{\ee}{\end{equation}}
\newcommand{\bea}{\begin{eqnarray}}
\newcommand{\eea}{\end{eqnarray}}
\newcommand{\bean}{\begin{eqnarray*}}
\newcommand{\eean}{\end{eqnarray*}}
\newcommand{\brray}{\begin{array}}
\newcommand{\erray}{\end{array}}
\newcommand{\ben}{\begin{equation}{nonumber}}
\newcommand{\een}{\end{equation}{nonumber}}

\newtheorem{dfn}{Definition}[section]
\newtheorem{thm}[dfn]{Theorem}
\newtheorem{lmma}[dfn]{Lemma}
\newtheorem{ppsn}[dfn]{Proposition}
\newtheorem{crlre}[dfn]{Corollary}
\newtheorem{xmpl}[dfn]{Example}
\newtheorem{rmrk}[dfn]{Remark}
\newcommand{\bdfn}{\begin{dfn}}
\newcommand{\bthm}{\begin{thm}}
\newcommand{\blmma}{\begin{lmma}}
\newcommand{\bppsn}{\begin{ppsn}}
\newcommand{\bcrlre}{\begin{crlre}}
\newcommand{\bxmpl}{\begin{xmpl}}
\newcommand{\brmrk}{\begin{rmrk}}
\newcommand{\edfn}{\end{dfn}}
\newcommand{\ethm}{\end{thm}}
\newcommand{\elmma}{\end{lmma}}
\newcommand{\eppsn}{\end{ppsn}}
\newcommand{\ecrlre}{\end{crlre}}
\newcommand{\exmpl}{\end{xmpl}}
\newcommand{\ermrk}{\end{rmrk}}

\newcommand{\IC}{\mathbb{C}}

\newcommand{\IN}{{I\! \! N}}

\newcommand{\IT}{\mathbb{T}}

\newcommand{\IZ}{\mathbb{Z}}

\newcommand{\cla}{{\cal A}}
\newcommand{\clb}{{\cal B}}
\newcommand{\clc}{{\cal C}}

\newcommand{\clg}{{\cal G}}
\newcommand{\clh}{{\cal H}}

\newcommand{\clk}{{\cal K}}
\newcommand{\cll}{{\cal L}}
\newcommand{\clm}{{\cal M}}

\newcommand{\clq}{{\cal Q}}

\newcommand{\cls}{{\cal S}}

\def\a*{{\cal A}_{h,*}}
\def\B{{\cal B}(h)}
\def\B1{{\cal B}_1(h)}
\def\b{{\cal B}^{\rm s.a.}(h)}
\def\b1{{\cal B}^{\rm s.a.}_1(h)}

\newcommand{\ot}{\otimes}

\newcommand{\raro}{\rightarrow}

\def \qed {$\Box$}

\begin{document}
	\[
\]
\begin{center}
{\large {\bf Some counterexamples in the theory of quantum isometry groups}}\\
by\\
{\large Jyotishman Bhowmick  }\\
and\\
{\large Debashish Goswami }\\

\end{center}

\begin{abstract}
    
  By considering spectral triples  on $S^{2}_{\mu, c}$ ( $c>0$) constructed by Chakraborty and Pal (\cite{chak_pal}), we show that 
    in general the quantum group of volume and orientation preserving isometries (in the sense of \cite{goswami2}) for a spectral triple of compact type  may not have a $C^*$-action,   and moreover, it can fail to be a matrix quantum group. It is also proved that  the category with objects consisting of those volume and orientation preserving quantum isometries which induce $C^*$-action on the $C^*$ algebra underlying the given spectral triple,  may not have a universal object.
  \end{abstract} 
  
 
\section{Introduction}
In a series of articles initiated by 
\cite{goswami} and followed by \cite{jyotish}, \cite{goswami2},  we have formulated and studied  a quantum group analogue of the group of Riemannian isometries of a classical or noncommutative manifold. This was motivated by  previous work of   a number of mathematicians including Wang, Banica, Bichon and others  (see, e.g. \cite{free}, \cite{wang}, 
\cite{ban1}, \cite{ban2}, \cite{bichon}, \cite{univ1} and references therein), who have defined quantum automorphism and quantum isometry groups  of finite spaces and finite dimensional algebras. Our theory of quantum isometry groups can be viewed as a natural generalization of such quantum automorphism or isometry groups of `finite' or `discrete' structures to the continuous or smooth set-up. Clearly, such a generalization is crucial to study the quantum symmetries in noncommutative geometry, and in particular, for a good understanding of quantum group equivariant spectral triples.

   The group of Riemannian isometries of a compact Riemannian manifold $M$ can be viewed as the universal object in the
      category of all compact metrizable groups acting on $M$, with smooth and isometric action. Moreover, assume that the manifold has a spin structure (hence in particular orientable, so we can fix a choice of orientation) and $D$ denotes the conventional Dirac operator acting as an unbounded self-adjoint operator on the Hilbert space $\clh$ of square integrable spinors. Then, it can be proved that the action of a compact group $G$ on the manifold lifts as a unitary representation (possibly of some group $\tilde{G}$ which is topologically a $2$-cover of  $G$, see \cite{CD} and \cite{Dabrowski_spinors} for more details) on the Hilbert space $\clh$ which commutes with $D$ if and only if   the action on the manifold is an orientation preserving isometric action. Therefore, to define the quantum analogue of the group of orientation-preserving Riemannian isometry group of a possibly noncommutative manifold given by a spectral triple $(\cla^\infty, \clh, D)$, 
       it is reasonable to  consider a
     category  ${\bf Q}^\prime(D)$ of compact quantum groups having unitary (co-) representation, say $U$, on $\clh$,  which commutes with $D$, and  the (a-priori von Neumann algebraic) action ${\rm ad}_U$  maps $\cla^\infty$ into its weak closure. 
A universal object in this category, if it exists, should define the `quantum group of orientation preserving Riemannian isometries' of the underlying spectral triple. Indeed (see \cite{goswami2}), if we consider a classical spectral triple, the subcategory of the category ${\bf Q}^\prime(D)$ consisting of groups has the classical group of orientation preserving isometries as the universal object, which justifies our definition of the quantum analogue.  Unfortunately, if we consider quantum group actions, even in the finite-dimensional (but with noncommutative $\cla$) situation  the category ${\bf Q}^\prime(D)$  may often fail to have a universal object.  It turns out, however, that if we fix any suitable faithful functional $ \tau_{R} $ on $\clb(\clh)$ (to be interpreted as the choice of a `volume form') then there exists a universal object in the subcategory $ {\bf Q}^{\prime}_{R}(D) $ of ${\bf Q}^\prime(D)$ obtained by restricting the object-class to the quantum group actions which also preserve the given functional.  The  subtle point to note here is that unlike the classical group actions on $\clb(\clh)$ which always preserve the usual trace, a quantum group action may not do so. In fact, it was proved by one of the authors in \cite{goswami_rmp} that given an object $(\clq, U)$ of ${\bf Q}^\prime(D)$ (where $\clq$ is the compact quantum group and $U$ denotes its unitary co-representation on $\clh$), we can find a suitable functional $\tau_{R}$ (which typically differs from the usual trace of $\clb(\clh)$ and can have a nontrivial modularity) which  is preserved by the action of $\clq$. This makes it quite natural to work in the setting  of twisted spectral data (as defined in \cite{goswami_rmp}). It may also be mentioned that in \cite{goswami2} we have actually worked in slightly bigger category $ {\bf Q}_{R}(D) $ of so called ` quantum family of orientation and volume preserving isometries '  and deduced that the universal object in $ {\bf Q}_{R}(D) $ exists and coincides with that of $ {\bf Q}^{\prime}_{R}(D) .$

Now, it is natural to ask whether ${\rm ad}_U$ induces a $C^*$-action (in the sense of Woronowicz and Podles) of the compact quantum group $QISO^+_R$ (or $QISO^+$ whenever it exists) on the $C^*$-completion $\cla$ of $\cla^\infty$. This is indeed the case for group actions on manifolds, where a measurable orientation preserving isometric action by a compact group  automatically becomes a smooth action, thanks to the Sobolev's theorem. However, we show in this article by giving a rather simple counterexample that such `automatic continuity' fails in the quantum situation: more precisely, ${\rm ad}_U$ does not in general induce a $C^*$-action of $QISO^+_R$ or $QISO^+$. From this, we also prove that the subcategory of ${\bf Q}^\prime_R$ (or ${\bf Q}^\prime$) consisting of  $(\tilde{Q}, U)$ for which ${\rm ad}_U$ induces $C^*$-action may fail to have a universal object, thus answering in the negative a question raised by S. Wang (private communication).

 \section{Notations and Preliminaries}
 \subsection{Basics of the theory of compact quantum groups}
 We begin by   recalling the definition of compact quantum groups and their actions from   \cite{woro}, \cite{woro1}.  A 
compact quantum group (to be abbreviated as CQG from now on)  is given by a pair $(\cls, \Delta)$, where $\cls$ is a unital separable $C^*$ algebra 
equipped
 with a unital $C^*$-homomorphism $\Delta : \cls \raro \cls \otimes \cls$ (where $\otimes$ denotes the injective tensor product)
  satisfying \\
  (ai) $(\Delta \ot id) \circ \Delta=(id \ot \Delta) \circ \Delta$ (co-associativity), and \\
  (aii) the linear spans of $ \Delta(\cls)(\cls \ot 1)$ and $\Delta(\cls)(1 \ot \cls)$ are norm-dense in $\cls \ot \cls$. \\
  It is well-known (see \cite{woro}, \cite{woro1}) that there is a canonical dense $\ast$-subalgebra $\cls_0$ of $\cls$, consisting of the matrix coefficients of
   the finite dimensional unitary (co)-representations (to be defined shortly) of $\cls$, and maps $\epsilon : \cls_0 \raro \IC$ (co-unit) and
   $\kappa : \cls_0 \raro \cls_0$ (antipode)  defined
    on $\cls_0$ which make $\cls_0$ a Hopf $\ast$-algebra. 

    A CQG $(\cls,\Delta)$ is said to (co)-act on a unital $C^*$ algebra $\clb$,
    if there is a  unital $C^*$-homomorphism (called an action) $\alpha : \clb \raro \clb \ot \cls$ satisfying the following :\\
    (bi) $(\alpha \ot id) \circ \alpha=(id \ot \Delta) \circ \alpha$, and \\
    (bii) the linear span of $\alpha(\clb)(1 \ot \cls)$ is norm-dense in $\clb \ot \cls$.\\
   
A unitary ( co ) representation of a CQG $ ( S, \Delta ) $ on a Hilbert space $ \clh $ is a map $ U $ from $ \clh $ to the Hilbert $\cls$ module $ \clh \otimes \cls $  such that the  element $ \widetilde{U} \in \clm ( \clk ( \clh ) \otimes \cls ) $ given by $\widetilde{U}( \xi \ot b)=U(\xi)(1 \ot b)$ ($\xi \in \clh, b \in \cls)$) is a unitary satisfying  $$ ({\rm  id} \otimes \Delta ) \widetilde{U} = {\widetilde{U}}_{(12)} {\widetilde{U}}_{(13)},$$ where for an operator $X \in \clb(\clh_1 \ot \clh_2)$ we have denoted by $X_{(12)}$ and $X_{(13)}$ the operators $X \ot I_{\clh_2} \in \clb(\clh_1 \ot \clh_2 \ot \clh_2)$, and $\Sigma_{23} X_{12} \Sigma_{23}$ respectively ($\Sigma_{23}$ being the unitary on $\clh_1 \ot \clh_2 \ot \clh_2$ which flips the two copies of $\clh_2$).
  
Given a unitary representation $U$ we shall denote by ${\rm ad}_U$ the $\ast$-homomorphism ${\rm ad}_U(X)=\widetilde{U}(X \ot 1){\widetilde{U}}^*$ for $X \in \clb(\clh)$.

We say that a (possibly unbounded) operator $T$ on $\clh$ commutes with $U$ if $T \ot I$ (with the natural domain) commutes with $\widetilde{U}$. Sometimes such an operator will be called $U$-equivariant.

 \subsection{The quantum group of orientation preserving Riemannian isometries}
 
 We briefly recall the definition of the quantum group of orientation preserving Riemannian isometries for a spectral triple (of compact type) $ ( \cla^{\infty}, \clh, D ) $ as in \cite{goswami2}.
 We consider the category ${\bf Q^\prime}(\cla^\infty, \clh, D) \equiv{\bf Q^\prime}(D)$ whose objects (to be called orientation preserving isometries) 
are the triplets $(\cls, \Delta, U)$, where    $(\cls, \Delta)$ is a CQG with a unitary representation $U$ in $\clh$, satisfying the following:\\
 (i) $U$ commutes with $D$, i.e. $\widetilde{U}$ commutes with $D \ot I$.\\ 
(ii) $({\rm id} \ot \phi) \circ {\rm ad}_U(a) \in ({\cla^\infty})^{\prime \prime}$ for all $ a \in \cla^\infty$ for every state $\phi$ on $\cls$, 
where ${\rm ad}_U(x):=\widetilde{U}(x \ot 1) {\widetilde{ U}}^* .  $

The category  $ {\bf {Q}^{\prime}}(D) $ may  not have a universal object in general, as seen in \cite{goswami2}. In case there is a universal object, 
we shall denote it by $ \widetilde{{QISO}^+}(D)$, with the corresponding representation $\widetilde{U}$, say, and we denote by $QISO^+(D)$ the Woronowicz 
subalgebra of $ \widetilde{{QISO}^+}(D) $ generated by the elements of the form $<\xi \ot 1, {\rm ad}_U(a)(\eta \ot 1)>$, 
where $\xi, \eta \in \clh, a \in \cla^\infty$ and $< \cdot, \cdot >$ is the $\widetilde{{QISO}^+}(D) $-valued inner product of $\clh \ot \widetilde{{QISO}^+}(D)$. The quantum group $QISO^+(D)$ will be called the quantum group of orientation-preserving Riemannian isometries of the spectral triple $(\cla^\infty, \clh, D)$. 

 Although the category ${\bf Q}^\prime(D)$ may sometimes fail to have a universal object,  we can always get a universal object in suitable subcategories which will 
be described now. Suppose that we are given an invertible positive (possibly unbounded) operator $R$ on $\clh$  which commutes with $D$.  Then we consider
 the full subcategory ${\bf Q}^\prime_R(D)$ of ${\bf Q}^\prime(D)$ by restricting the object class to those $(\cls, \Delta, U)$ for which $\cls$ is $\tau_R$-invariant, i.e.  $(\tau_R \ot {\rm id})({\rm ad}_U(X))=\tau_R(X)1$ for all $X$ in the $\ast$-subalgebra generated by operators of the form $|\xi><\eta|$, 
where $\xi, \eta$ are eigenvectors of the operator $D$ which by assumption has discrete spectrum, and $\tau_R(X)={\rm Tr}(RX)=<\eta, R \xi>$ 
for $X=|\xi><\eta|$.  We shall call the objects of ${\bf Q}^\prime_R (D)$ orientation and ($R$-twisted) volume preserving isometries.
It is clear (see Remark 2.9 in \cite{goswami2}) that when $Re^{-tD^2}$ is trace-class for some $t>0$, the above condition is equivalent to the condition that ${\rm ad}_U$ preserves the bounded normal functional ${\rm Tr}(\cdot ~Re^{-tD^2})$ on the whole of $\clb(\clh)$.  It is shown in \cite{goswami2} that the category ${\bf Q}^\prime_R(D)$ always admits a universal object,  to be denoted by $ \widetilde{{QISO}^+_R}(D)$, and the Woronowicz subalgebra generated  by $\{ <\xi \ot 1, {\rm ad}_W(a) (\eta \ot 1)>,~\xi, \eta \in \clh, a \in \cla^\infty \}$ (where $W$ is the unitary representation of $\widetilde{{QISO}^+_R}(D)$ in $\clh$) will be denoted by $QISO^+_R(D)$ and called the quantum group of orientation and ($R$-twisted) volume preserving Riemannian isometries of the spectral triple. 

It may also be noted that, if $QISO^+(D)$ exists, it will also be isomorphic with $QISO^+_R(D)$ for some (not necessarily unique) $R$. In fact, $QISO^+(D) \cong QISO^+_R(D)$ for any $R$ such that $\widetilde{QISO^+}(D)$ is $\tau_R$-invariant.

\section{Quantum isometry group of the spectral triple of Chakraborty and Pal on $ S^{2}_{\mu,c}, c > 0 $}

The Podles' sphere $ S^{2}_{\mu,c} $ is the  universal $ C^{*} $ algebra generated by elements $ A $ and $ B $ satisfying the relations:
$$ A^{*} = A,~  AB = \mu^{-2} BA,$$
$$ B^* B = A - A^2 + cI, ~ B B^* = \mu^2 A - \mu^4 A^2 + cI. $$


  Let us describe the spectral triple on $S^2_{\mu,c}$  introduced and studied  in \cite{chak_pal}.
  
Let $ \clh_{+} = \clh_{-} = l^2 ( \IN \bigcup\{0\} ),\clh = \clh_{+} \oplus \clh_{-} .$

Let $\{  e_n, n \geq 0 \} $ be the canonical  orthonormal basis of $ \clh_{+} = \clh_{-} $ and $ N $ be the operator defined on it by $ N ( e_n ) = n e_n .$

We recall the irreducible  representations $ \pi_{\pm}  : \clh_{\pm} \rightarrow \clh_{\pm} $  given by:
$$ \pi_{\pm} ( A ) e_n = \lambda_{\pm} {\mu}^{2n} e_n, $$
$$ \pi_{\pm} ( B ) e_n = {c_{\pm}( n )}^{\frac{1}{2}} e_{n - 1},  $$
where $ e_{- 1} = 0,~ \lambda_{\pm} =  \frac{1}{2} \pm {( c + \frac{1}{4} )}^{\frac{1}{2}}, $
\be \label{sphere_chak_pal_c_pm(n)} c_{\pm} ( n ) = \lambda_{\pm} \mu^{2n} - {( \lambda_{\pm} \mu^{2n} )}^{2} + c. \ee
Let $ \pi = \pi_{+} \oplus \pi_{-} $ and $ D = \left ( \begin {array} {cccc}
   0 & N  \\ N & 0 \end {array} \right ) .$
   
Then $ ( S^{2}_{\mu,c}, \pi, \clh, D ) $ is a spectral triple of compact type.

We note that the eigenvalues of $ D $ are $ \{ m : m \in \IZ \} $ and the  one-dimensional eigenspace corresponding to $m$ is spanned by $ \left ( \begin {array} {cccc}
   e_m \\ e_m \end {array} \right ) $ when $m \geq 0$,  and by $ \left ( \begin {array} {cccc}
   e_{-m} \\ - e_{-m} \end {array} \right ) $ when $m<0$.
   
 It follows from the definition of $ \pi_{\pm} ( B ) $ that 
  $$ \pi_{+} ( B^* ) ( e_n ) = {c_{+} ( n + 1 )}^{\frac{1}{2}} e_{ n + 1 }, $$  
  $$ \pi_{-} ( B^* ) ( e_n ) = {c_{-} ( n + 1 )}^{\frac{1}{2}} e_{ n + 1 }. $$

\blmma

\label{sphere_chak_pal_algebra_preservation}
 
$ \pi {( S^{2}_{\mu,c} )}^{''} = \{ \left ( \begin {array} {cccc}
   X_{11} & X_{12}  \\ X_{21} & X_{22} \end {array} \right ) \in \clb ( \clh \oplus \clh ) : X_{12} = X_{21} = 0 \}. $
   
 \elmma
 
 {\it Proof :} It suffices to prove that the commutant $ {\pi ( S^{2}_{\mu,c}} )^{\prime} $ is the von Neumann algebra of operators of the form $  \left ( \begin {array} {cccc}
   c_1 I & 0  \\ 0 & c_2 I \end {array} \right ) $ for  $ c_1, c_2 \in \IC.$ We use the fact that $ \pi_{+} $ and $ \pi_{-} $ are irreducible representations.
 
      Let $ X = \left ( \begin {array} {cccc}
   X_{11} & X_{12}  \\ X_{21} & X_{22} \end {array} \right ) \in { \pi ( S^{2}_{\mu,c} )}^{\prime} .$ Using the fact that $ X $ commutes with $ \pi ( A ),~ \pi ( B ),~ \pi ( B^* ) ,$ we have: $ X_{11} \in { \pi_{+} ( S^{2}_{\mu,c} )}^{\prime} \cong \IC $ and $X_{22} \in {\pi_{-} (S^{2}_{\mu,c})}^{\prime} \cong \IC ,$ so let $ X_{11} = c_1 I, ~ X_{22} = c_2 I $ for some $ c_1,c_2 .$
   
   Moreover, \be \label{sphere_chak_pal_i} X_{12} \pi_{-} ( A ) = \pi_{+} ( A ) X_{12}, \ee   
             \be \label{sphere_chak_pal_ii} X_{12} \pi_{-} ( B ) = \pi_{+} ( B ) X_{12}, \ee             
             \be \label{sphere_chak_pal_iii} X_{12} \pi_{-} ( B^* ) = \pi_{+} ( B^* ) X_{12}. \ee             
   Now, ( \ref{sphere_chak_pal_ii} ) implies $ X_{12} e_0 \in Ker ( \pi_{+} ( B ) ) = \IC e_0.$
   Let $ X_{12} e_0= p_0 e_0.$ 
   
   We have, $ \pi_{+} ( B ) ( X_{12} e_1 ) = {c_{-}( 1)}^{\frac{1}{2}} X_{12} e_0 = p_0 {c_{-}( 1 )}^{\frac{1}{2}} e_0 ,$ i.e, $ \pi_{+} ( B ) ( X_{12} e_1 ) \in \IC e_0.$
   
 Since it follows from the definition of $ \pi_{+} ( B ) $ that $ \pi_{+} ( B ) $ maps $ \overline{{\rm span} ~ \{ e_i : i \geq 2 \}} $ to $ {( \IC e_0 )}^{\bot} = \overline{{\rm span} \{ e_i : i \geq 1 \}} , ~ X_{12} e_1 $ must belong to $ span \{ e_0, ~ e_1 \}. $
   
   Inductively, we conclude that  for all $ n,~ X_{12} ( e_n ) \in {\rm span} \{ e_0, e_1,......e_n \}.$
   
   Using the definition of $ \pi_{\pm} ( B^* ) e_n $ along with ( \ref{sphere_chak_pal_iii} ), we have  $ {c_{-} ( 1 )}^{\frac{1}{2}} X_{12} e_1 = p_0 {c_{+} ( 1 )}^{\frac{1}{2}} e_1 ,$ i.e, $ X_{12} e_1 = p_0 \frac{{c_{+} ( 1)}^{\frac{1}{2}}}{{c_{-} ( 1 )}^{\frac{1}{2}}} e_1 .$ 
   
   We argue in a similar way  by induction that $ X_{12} e_n = c^{\prime}_{n} e_n $ for some constants $ c^{\prime}_{n}.$
   
   Now we apply  ( \ref{sphere_chak_pal_iii} ) and ( \ref{sphere_chak_pal_ii} ) on the vectors $ e_n$ and $e_{n + 1}$ to get 
   
  $ c^{\prime}_{n + 1} = \frac{c^{\prime}_{n} {c_{+} ( n + 1 )}^{\frac{1}{2}}}{{c_{-} ( n + 1 )}^{\frac{1}{2}}} $ and $ c^{\prime}_{n + 1} = \frac{c^{\prime}_{n} {c_{-} ( n + 1 )}^{\frac{1}{2}}}{{c_{+} ( n + 1 )}^{\frac{1}{2}}}. $
  
 Since $ c_{+} ( n + 1 ) \neq c_{-} ( n + 1 ) $ for $ n \geq 0,$ we have $ c^{\prime}_{n} = 0 .$
  
  Hence, $ c^{\prime}_{n}  = 0 $ for all $ n $ implying $ X_{12} = 0.$
  
  It follows similarly that $ X_{21} = 0 .$ \qed
   
   \vspace{4mm}        
 
 Let $ ( \widetilde{\clq}, \Delta, U ) $ be an object in the category $ {\bf Q^\prime}(D) ,$ with  $ \alpha={\rm ad}_U$ and  the corresponding Woronowicz $ C^{*} $ subalgebra of $ \widetilde{\clq} $ generated by $\{ <(\xi \ot 1), \alpha(x)(\eta \ot 1)>_{\tilde{\clq}},~\xi, \eta \in \clh,~x \in S^2_{\mu,c} \}$ is denoted by $\clq$. Assume, without loss of generality, that the representation $U$ is faithful. 
 
 Since $ U $ commutes with $ D $, it preserves the eigenvectors  $ \left ( \begin {array} {cccc}
   e_n \\ e_n \end {array} \right ) $ and  $ \left ( \begin {array} {cccc}
   e_n \\ - e_n \end {array} \right ) .$
   
   Let $ U ~  \left ( \begin {array} {cccc}
   e_n \\ e_n \end {array} \right ) ~  = ~  \left ( \begin {array} {cccc}
   e_n \\ e_n \end {array} \right )  \otimes q^{+}_{n}, $
   
   $ U ~  \left ( \begin {array} {cccc}
   e_n \\ - e_n \end {array} \right ) ~  = ~ \left ( \begin {array} {cccc}
   e_n \\ - e_n \end {array} \right )  \otimes q^{-}_{n}, $
   
   for some $ q^{+}_{n}, q^{-}_{n} \in \widetilde{\clq} .$ By the assumption of faithfulness, the unital $C^*$ algebra generated by the elements $q^+_n, q^-_n$ will be the whole of $\tilde{\clq}$.

 \blmma
 
 \label{sphere_chak_pal_alg}
 We have:
\be \label{sphere_chak_pal_alg1} q^{+}_{n}q^{-*}_{n} = q^{-}_{n}q^{+*}_{n} ~ {\rm  for ~ all} ~ n, \ee

\be \label{sphere_chak_pal_alg2} ( {c_{+}( n )}^{\frac{1}{2}} + {c_{-}( n )}^{\frac{1}{2}} ) ( q^{+}_{n - 1} q^{+*}_{n} - q^{-}_{n - 1} q^{-*}_{n} ) +  ( {c_{+}( n )}^{\frac{1}{2}} - {c_{-}( n )}^{\frac{1}{2}} ) ( q^{+}_{n - 1} q^{-*}_{n} - q^{-}_{n - 1} q^{+*}_{n} ) = 0 ~ {\rm  for ~ all} ~ n \geq 1, \ee

\be \label{sphere_chak_pal_alg3} ( {c_{+}( n )}^{\frac{1}{2}} + {c_{-}( n )}^{\frac{1}{2}} ) ( q^{+}_{n - 1} q^{+*}_{n} - q^{-}_{n - 1} q^{-*}_{n} ) +  ( {c_{+}( n )}^{\frac{1}{2}} - {c_{-}( n )}^{\frac{1}{2}} ) ( q^{-}_{n - 1} q^{+*}_{n} - q^{+}_{n - 1} q^{-*}_{n} ) = 0 ~ {\rm  for ~ all} ~ n \geq 1, \ee

\be \label{sphere_chak_pal_alg4} ( {c_{+}( n + 1 )}^{\frac{1}{2}} + {c_{-}( n + 1 )}^{\frac{1}{2}} ) ( q^{+}_{n + 1} q^{+*}_{n} - q^{-}_{n + 1} q^{-*}_{n} ) =  ( {c_{+}( n + 1 )}^{\frac{1}{2}} - {c_{-}( n + 1 )}^{\frac{1}{2}} ) ( q^{-}_{n + 1} q^{+*}_{n} - q^{+}_{n + 1} q^{-*}_{n} )  ~ {\rm  for ~ all}~ n.  \ee
 
 \elmma
 
{\it Proof :} One has,
 $ \alpha ( A )  \left ( \begin {array} {cccc}
   e_n \\ 0 \end {array} \right ) $
   
   $ = \widetilde{U} ( A \otimes 1 ) {\widetilde{U}}^{*} ( \left ( \begin {array} {cccc}
   e_n \\ 0 \end {array} \right ) \otimes 1 ) $
   
   $ = \frac{1}{2} \widetilde{U} ( A \otimes 1 ) {\widetilde{U}}^{*} ( \left ( \begin {array} {cccc}
   e_n \\ e_n \end {array} \right ) \otimes 1 + \left ( \begin {array} {cccc}
   e_n \\ - e_n \end {array} \right ) \otimes 1 ) $
   
   $ = \frac{1}{2} \widetilde{U} [ \left ( \begin {array} {cccc}
  \lambda_{+} \mu^{2n} e_n \\ \lambda_{-} \mu^{2n} e_n \end {array} \right ) \otimes q^{+*}_{n} + \left ( \begin {array} {cccc}
   \lambda_{+} \mu^{2n}e_n \\ - \lambda_{-} \mu^{2n} e_n \end {array} \right ) \otimes q^{-*}_{n} ] $
   
   $ = \frac{1}{4} \widetilde{U} [ \left ( \begin {array} {cccc}
   e_n \\  e_n \end {array} \right ) \otimes \{ \lambda_{+} \mu^{2n} ( q^{+*}_{n} + q^{-*}_{n} ) + \lambda_{-} \mu^{2n} ( q^{+*}_{n} - q^{-*}_{n} ) \} + \left ( \begin {array} {cccc}
   e_n \\ - e_n \end {array} \right ) \otimes \{ \lambda_{+} \mu^{2n} ( q^{+*}_{n} + q^{-*}_{n} ) - \lambda_{-} \mu^{2n} ( q^{+*}_{n} - q^{-*}_{n} ) \}  ] $
   
   $ =  \left ( \begin {array} {cccc}
   e_n \\ 0 \end {array} \right ) \otimes \frac{1}{4} \{ \lambda_{+} \mu^{2n} ( 1 + q^{+}_{n} q^{-*}_{n} ) + \lambda_{-} \mu^{2n} ( 1 - q^{+}_{n} q^{-*}_{n} ) + \lambda_{+} \mu^{2n}( 1 + q^{-}_{n} q^{+*}_{n} ) - \lambda_{-}\mu^{2n} (  q^{-}_{n} q^{+*}_{n} - 1 ) \} + \left ( \begin {array} {cccc}
   0 \\ e_n \end {array} \right ) \otimes \frac{1}{4} \{ \lambda_{+} \mu^{2n} ( 1 + q^{+}_{n} q^{-*}_{n} ) + \lambda_{-} \mu^{2n} ( 1 - q^{+}_{n} q^{-*}_{n} ) - \lambda_{+} \mu^{2n} ( 1 + q^{-}_{n} q^{+*}_{n} ) + \lambda_{-}\mu^{2n} (  q^{-}_{n} q^{+*}_{n} - 1 ) \}.  $

 Since $ \alpha ( A ) $ maps $ \pi( S^{2}_{\mu,c} ) $ into its double commutant, we conclude by using the description of $ {\pi ( S^{2}_{\mu,c} )}^{\prime \prime} $ given in Lemma \ref{sphere_chak_pal_algebra_preservation} that the coefficient of $ \left ( \begin {array} {cccc}
   0 \\ e_n \end {array} \right ) $ in $ \alpha ( A ) \left ( \begin {array} {cccc}
   e_n \\ 0 \end {array} \right ) $ must be $ 0 ,$ which implies $ q^{+}_{n}q^{-*}_{n} = q^{-}_{n}q^{+*}_{n} .$
   
 Proceeding in a similar way, ( \ref{sphere_chak_pal_alg2} ),( \ref{sphere_chak_pal_alg3} ), ( \ref{sphere_chak_pal_alg4} ) follow  from the facts that coefficients of $  \left ( \begin {array} {cccc}
   0 \\ e_{n - 1} \end {array} \right ) , \left ( \begin {array} {cccc}
   e_{n - 1} \\ 0 \end {array} \right ) $ and $ \left ( \begin {array} {cccc}
   0 \\ e_{n + 1} \end {array} \right ) $ in $ \alpha ( B )  \left ( \begin {array} {cccc}
   e_n \\ 0 \end {array} \right ) ,$   
  $  \alpha ( B ) \left ( \begin {array} {cccc}
   0 \\ e_{n} \end {array} \right ) $ and   $ \alpha ( B^* ) \left ( \begin {array} {cccc}
   e_n \\ 0 \end {array} \right ) $ ( respectively ) are zero. \qed
\vspace{3mm}\\
For a series of the form $\sum_n T_n,$ with $T_n \in \cll(E)$ for some Hilbert $C^*$-module $E$, let us say that $\sum_n T_n$ converges point-wise if  for every $\xi$ in $E$, the series $\sum_n T_n \xi$ converges in the norm of $E$. With this terminology, we have the following 
   
   \bcrlre
   
   \label{sphere_chak_pal_corollary}
   
   If $ P_n, Q_n $ denote the projections onto the subspace generated by $ \left ( \begin {array} {cccc}
   e_n \\ 0 \end {array} \right ) $ and $ \left ( \begin {array} {cccc}
   0 \\ e_n \end {array} \right ) $ respectively, then $P_n, Q_n \in \pi(S^2_{\mu,c})$ and we have 
   
   \bean \alpha ( A ) = \sum^{\infty}_{n = 0} A P_n \otimes \frac{1}{2 \lambda_{+}} \{ \lambda_{+} ( 1 + q^{+}_{n} q^{-*}_{n} ) + \lambda_{-} ( 1 - q^{+}_{n} q^{-*}_{n} ) \} \eean
   
   \bean + \sum^{\infty}_{n = 0} A Q_n \otimes \frac{1}{2 \lambda_{-}} \{ \lambda_{+} ( 1 - q^{+}_{n} q^{-*}_{n} ) + \lambda_{-} ( 1 + q^{+}_{n} q^{-*}_{n} ) \}, \eean
     
    \bean \alpha ( B ) = \sum^{\infty}_{n = 1} B P_n \otimes \frac{1}{4 c_{+} ( n ) } [ ( {c_{+} ( n )}^{\frac{1}{2}} + {c_{-} ( n )}^{\frac{1}{2}} ) ( q^{+}_{n - 1} q^{+*}_{n} + q^{-}_{n - 1} q^{-*}_{n} ) \eean
    
    \bean + ( {c_{+} ( n )}^{\frac{1}{2}} - {c_{-} ( n )}^{\frac{1}{2}} ) ( q^{-}_{n - 1} q^{+*}_{n} + q^{+}_{n - 1} q^{-*}_{n} ) ] + \sum^{\infty}_{n = 1} B Q_n \otimes \frac{1}{4 c_{-} ( n ) } [ ( {c_{+} ( n )}^{\frac{1}{2}} + {c_{-} ( n )}^{\frac{1}{2}}  ) \eean
    
    \bean ( q^{+}_{n - 1} q^{+*}_{n} + q^{-}_{n - 1} q^{-*}_{n} ) - ( {c_{+} ( n )}^{\frac{1}{2}} - {c_{-} ( n )}^{\frac{1}{2}} ) ( q^{+}_{n - 1} q^{-*}_{n} + q^{-}_{n - 1} q^{+*}_{n} ) ]. \eean
    
    Here, each of the series of elements of $\cll(\clh \ot \tilde{\clq})$ appearing on the right hand sides,  converges in the point-wise sense defined before the statement of the corollary.
   
   \ecrlre
   
   {\it Proof :} It is easy to see that $P_n$, $Q_n$ are the eigenprojections of $B^*B$ for the eigenvalues  $ c_{+} ( n ) , c_{-}( n )  $ respectively, and these are isolated points in the spectrum, that is in the  point spectrum, so by the spectral theorem, $P_n, Q_n$ belong to $C^*(B^*B) \subset \pi(S^2_{\mu,c}).$ 
The proof now follows by applying Lemma \ref{sphere_chak_pal_alg} on $ \alpha ( A ) \left ( \begin {array} {cccc}
   e_n \\ 0 \end {array} \right ) , \alpha ( A ) \left ( \begin {array} {cccc}
   0 \\ e_{n} \end {array} \right ) ,~  \alpha ( B ) \left ( \begin {array} {cccc}
   e_n \\ 0 \end {array} \right ) , \alpha ( B ) \left ( \begin {array} {cccc}
   0 \\ e_{n} \end {array} \right ).$  \qed

 \bppsn
   
   \label{sphere_chak_pal_description of Q'+}

   As a $C^*$-algebra $ \widetilde{\clq}$ is  generated by the unitaries $  \{ q^{+}_{n} \}_{n \geq 0}$, and the self-adjoint unitary $y_0  = q^{-*}_{0} q^{+}_{0} .$ Moreover, $ \clq $ is generated by unitaries $ z_n = q^{+}_{n - 1} q^{+*}_{n}, ~ n \geq 1 $ along with a self adjoint unitary $ w^{\prime} .$ 
   
   \eppsn

 {\it Proof :} Replacing n + 1 by n in ( \ref{sphere_chak_pal_alg4} ) we have, 
\be  \label{sphere_chak_pal_id5} ( {c_{+}( n  )}^{\frac{1}{2}} + {c_{-}( n )}^{\frac{1}{2}} ) ( q^{+}_{n} q^{+*}_{n - 1} - q^{-}_{n} q^{-*}_{n- 1} ) - ( {c_{+}( n )}^{\frac{1}{2}} - {c_{-}( n )}^{\frac{1}{2}} ) ( q^{-}_{n} q^{+*}_{n - 1} - q^{+}_{n} q^{-*}_{n - 1} ) = 0 ~ {\rm for ~ all } ~ n \geq 1. \ee
Subtracting ( \ref{sphere_chak_pal_id5} ) from the equation obtained by applying $ \ast  $ on ( \ref{sphere_chak_pal_alg2} ), we have : $ 2 ( {c_{+}( n  )}^{\frac{1}{2}} - {c_{-}( n )}^{\frac{1}{2}} ) ( q^{-}_{n} q^{+*}_{n - 1} - q^{+}_{n} q^{-*}_{n - 1} ) = 0 $  for all $ n \geq 1  $ which gives :
\be \label{sphere_chak_pal_id6} q^{-}_{n} q^{+*}_{n - 1} = q^{+}_{n} q^{-*}_{n - 1}   ~ {\rm  for ~ all} ~  n \geq 1. \ee
Using ( \ref{sphere_chak_pal_id6} ) in ( \ref{sphere_chak_pal_id5} ), we have
 \be \label{sphere_chak_pal_id7} q^{+}_{n} q^{+*}_{n - 1} = q^{-}_{n} q^{-*}_{n - 1} ~ {\rm  for ~ all} ~ n \geq 1. \ee 
 Let $ y_n = q^{-*}_{n} q^{+}_{n}.$
 
 Then, using  ( \ref{sphere_chak_pal_alg1} ),  we observe that $ y_n $ is a self adjoint unitary.
 
 Moreover, from ( \ref{sphere_chak_pal_id6} ), we have 
  \be \label{sphere_chak_pal_id8} q^{-}_{n} = q^{+}_{n} y_{n - 1} ~ {\rm  for ~ all} ~ n \geq 1. \ee  
  From ( \ref{sphere_chak_pal_alg4} ) we have  
  \be \label{sphere_chak_pal_alg4a} ( {c_{+}( n )}^{\frac{1}{2}} + {c_{-}( n )}^{\frac{1}{2}} ) ( q^{+}_{n } q^{+*}_{n - 1} - q^{-}_{n} q^{-*}_{n - 1} ) =  ( {c_{+}( n )}^{\frac{1}{2}} - {c_{-}( n )}^{\frac{1}{2}} ) ( q^{-}_{n} q^{+*}_{n - 1} - q^{+}_{n} q^{-*}_{n - 1} )  ~ {\rm  for ~ all} ~ n \geq 1.  \ee  
Next, using ( \ref{sphere_chak_pal_id7} ), we obtain  $ q^{-*}_{n} q^{+}_{n} = q^{-*}_{n - 1} q^{+}_{n - 1} $ for all $n \geq 1$ implying  
  \be \label{sphere_chak_pal_id9} y_n = y_{n - 1} ~ {\rm  for ~ all} ~ n \geq 1. \ee  
  From ( \ref{sphere_chak_pal_id8} ) and ( \ref{sphere_chak_pal_id9} ) and the faithfulness of the representation $U$, we conclude that $ \widetilde{\clq} $ is generated by $ \{ q^{+}_{n} \}_{n \geq 0} $ and $ y_{0}.$

 Now we prove the second part of the proposition.
 Let $ \widetilde{P_n} = P_n + Q_n \in \pi(S^2_{\mu, c}).$
  Then, for all $ v \in \clh, \alpha (  \widetilde{P_n} )v =  \widetilde{P_n}v \otimes 1.$

   We note that for all $ v \in \clh, ~ \alpha ( A \widetilde{P_k} )v = \alpha ( A ) ( \widetilde{P_k}v \otimes 1 ) = A P_k v \otimes \frac{1}{2 \lambda_{+}} \{ \lambda_{+} ( 1 + q^{+}_{k}q^{-*}_{k} ) + \lambda_{-} ( 1 - q^{+}_{k} q^{-*}_{k} ) \} +  A Q_k v \otimes \frac{1}{2 \lambda_{-}} \{ \lambda_{+} ( 1 - q^{+}_{k} q^{-*}_{k} ) + \lambda_{-} ( 1 + q^{+}_{k} q^{-*}_{k} ) \} .$
   
   Hence, $ \lambda_{+} ( 1 + q^{+}_{m}q^{-*}_{m} ) + \lambda_{-} ( 1 - q^{+}_{m} q^{-*}_{m} ) $ and $ \lambda_{+} ( 1 - q^{+}_{m} q^{-*}_{m} ) + \lambda_{-} ( 1 + q^{+}_{m} q^{-*}_{m} ) \in \clq $  for all $ m.$
   
   Subtracting, we get $ q^{+}_{m} q^{-*}_{m} \in \clq.$ 
   
 Using the expression of $ \alpha ( B ) $ in a similar way, we have 
   
   \be \label{sphere_chak_pal_id10}   q^{+}_{n - 1} q^{+*}_{n} + q^{-}_{n - 1} q^{-*}_{n} \in \clq \mbox{ for  all} ~  n \geq 1,  \ee    
   \be \label{sphere_chak_pal_id11}  q^{-}_{n - 1} q^{+*}_{n} + q^{+}_{n - 1} q^{-*}_{n} \in \clq \mbox{ for  all} ~  n \geq 1.  \ee   
  Using $ q^{-}_{n} = q^{+}_{n} y_0 $ ( as observed from the first part of the proof ) in ( \ref{sphere_chak_pal_id10} ) and ( \ref{sphere_chak_pal_id11} ), we conclude that $ q^{+}_{n - 1} q^{+*}_{n} $ and $ q^{+}_{n - 1} y_0 q^{+*}_{n}  \in \clq$ for all $ n \geq 1. $
  
  Let $$ z_n = q^{+}_{n - 1} q^{+*}_{n}, $$  
      $$ w_n =  q^{+}_{n - 1} y_0 q^{+*}_{n}, $$      
       for all $ n \geq 1. $
      
      Then, we observe that $ z^{*}_{n} w_n = q^{+}_{n} y_0 q^{+*}_{n} = q^{+}_{n} q^{-*}_{n}. $
      
      Moreover, \bean q^{+}_0 q^{-*}_0 = q^{+}_0 {( q^{+}_0 y^{*}_0 )}^{*} = q^{+}_0 y_0 q^{+*}_0 = q^{+}_0 y_0 q^{+*}_1  q^{+}_1 q^{+*}_0 = w_1 z^{*}_1.\eean      
    Thus  for all $ n \geq 0, ~ q^{+}_{n} q^{-*}_{n} $ belong to $ C^{*} ( \{ z_n, w_n \}_{n \geq 1} ) .$     
   \bean {\rm~ For ~ all} ~ n \geq 2, ~ q^{-}_{n - 1} q^{-*}_n = q^{+}_{n - 1} y^{*}_{n - 2} {( q^{+}_n y^{*}_{n -1} )}^* = q^{+}_{n - 1} y^*_0 y_0 q^{+*}_n = q^{+}_{n - 1} q^{+*}_n = q^{+}_{n - 1} q^{+*}_n = z_n, \eean   
   \bean q^{-}_0 q^{-*}_1 = q^{+}_0 y_0 {( q^{+}_1 y_0)}^{*} = q^{+}_0 y_0 y^{*}_0 q^{+*}_1 = q^{+}_0 q^{+*}_1 = z_1.\eean
   
   Finally,\bean q^{-}_{n - 1} q^{+*}_n = q^{+}_{n - 1} y^*_0 q^{+*}_n = w_n \eean  and \bean q^{-}_0 q^{+*}_1 = q^{+}_0 y_0 q^{+*}_1 = w_1. \eean

    Now, from the expressions of $ \alpha ( A ) $ and $\alpha( B ) ,$ it is clear that $ \clq $ is generated by $ q^{+}_{n} {q^{-}_{n}}^{*}, ~ q^{+}_{n - 1} {q^{+}_{n}}^{*} + q^{-}_{n - 1} {q^{-}_{n}}^{*}, ~ q^{-}_{n - 1} {q^{+}_{n}}^{*} + q^{+}_{n - 1} {q^{-}_{n}}^{*} .$ By the above made observations, these belong to $ C^{*} ( {\{ z_n, w_n \}}_{n \geq 1} ) $ which implies that $ \clq $ is a $ C^* $ subalgebra of $ C^{*} ( {\{ z_n, w_n \}}_{n \geq 1} ) .$ Moreover, from the definitions of $ z_n,~ w_n $ it is clear that $ C^{*} ( {\{ z_n, w_n \}}_{n \geq 1} ) $ is a $ C^* $ subalgebra of $ \clq .$ 
    
   Therefore, $ \clq \cong C^{*} ( \{ z_n, w_n \}_{n \geq 1} ) .$
    
    In fact, a simpler description is possible by noting that $  w_{n + 1} = {z_n}^{*} w_n z_{n + 1} $ which implies  $ \{ w_n \}_{n \geq 1} \in C^{*} ( \{ z_n \}_{n \geq 1}, ~ w_1 ) .$
    
    Defining $ w^{\prime} = w^{*}_{1} z_1 $ and using the expressions for $ z_1 $ and $ w_1 $ we note that $ w^{\prime} $ is a self adjoint unitary. 
    
    Thus $ \clq \cong C^{*} \{ \{ z_n \}_{n \geq 1}, ~ w^{\prime} \}.$ \qed
    
    \vspace{4mm}


 We now look at the coproduct, say $\Delta$ of $\tilde{\clq}.$ We have,    
     $$ ( id \otimes \Delta ) U  \left ( \begin {array} {cccc}
   e_n \\ e_n \end {array} \right ) = U_{(12)} U_{(13)}  \left ( \begin {array} {cccc}
   e_n \\ e_n \end {array} \right ), $$
   from which we immediately conclude the following:  
    \blmma
    
    \label{sphere_chak_pal_coproduct}
    
   $ \Delta ( q^{\pm}_{n} ) = q^{\pm}_{n} \otimes q^{\pm}_{n}, $
    
    $ \Delta ( y_1 ) = y_1 \otimes y_1. $
    
    \elmma

    
   
   
   
   

 Let us now consider the quantum group $\tilde{\clq}_0$ which is the dual of the discrete group $\IZ_2 \ast \IZ^\infty$, where $\IZ^\infty=\IZ \ast \IZ \ast \cdots$ denotes the free product of countably infinitely many copies of $\IZ$. The underlying $C^*$-algebra is nothing but $C(\IZ_2) \ast C(\IT) \ast C(\IT) \ast  \cdots$, and let us denote by  $ r^{+}_{n} $ the generator of n th copy of $ C ( \IT ) $ and by $ y $  the generator of $ C ( \IZ_{2} ).$

Define 
 \bean V \left ( \begin {array} {cccc}
   e_n \\ e_n \end {array} \right ) = \left ( \begin {array} {cccc}
   e_n \\ e_n \end {array} \right ) \otimes r^{+}_{n}, \eean   
 \bean  V \left ( \begin {array} {cccc}
   e_n \\ - e_n \end {array} \right ) = \left ( \begin {array} {cccc}
   e_n \\ - e_n \end {array} \right ) \otimes r^{+}_{n} y. \eean
   Then $ V $ commutes with $ D $ and it can be easily checked that $ V $ is a unitary representation of $\tilde{\clq}_0$ , i.e. $(\tilde{\clq}_0, \Delta_0, V)$ (where $\Delta_0$ denotes the coproduct on $\tilde{\clq}_0$, given by $\Delta_0(r^+_n) = r^+_n \ot r^+_n,~\Delta_0(y)=y \ot y$) is an object in ${\bf Q}^\prime(D).$
   
   Setting $ r^{-}_{n} = r^{+}_{n} y ,$ we observe that 
    $ r^{-}_{n} $ is a unitary and satisfies :   
   \be \label{sphere_chak_pal_id1} r^{-}_{n} r^{+*}_{n - 1} = r^{+}_{n} r^{-*}_{n - 1} ~~{\rm  for ~ all} ~ n \geq 1, \ee   
   \be \label{sphere_chak_pal_id2} r^{+}_{n} r^{+*}_{n - 1} = r^{-}_{n} r^{-*}_{n - 1} ~~ {\rm  for ~ all} ~ n \geq 1.   \ee   
   Using ( \ref{sphere_chak_pal_id1} ) in ( \ref{sphere_chak_pal_id2} ) we have $ r^{+}_{n - 1} = r^{-}_{n - 1} r^{+*}_{n - 1} r^{-}_{n - 1} .$
   
   This implies   
   \be \label{sphere_chak_pal_r1} r^{+}_{n}{r^{-}_{n}}^{*} = r^{-}_{n} r^{+*}_{n} ~~ {\rm  for ~ all} ~ n. \ee    
     Moreover, taking $ \ast $ on ( \ref{sphere_chak_pal_id1} ) and ( \ref{sphere_chak_pal_id2} )  respectively, we get the following:     
    \be \label{sphere_chak_pal_r2} r^{+}_{n - 1} r^{-*}_{n} - r^{-}_{n - 1}r^{+*}_{n} = 0 ~~ {\rm  for ~ all} ~ n \geq 1, \ee    
    \be   \label{sphere_chak_pal_r3}  r^{+}_{n - 1} r^{+*}_{n} - r^{-}_{n - 1} r^{-*}_{n} = 0 ~~ {\rm  for ~ all} ~ n \geq 1. \ee    
    Thus, the equations ( \ref{sphere_chak_pal_alg1} ) -  ( \ref{sphere_chak_pal_alg4} ) in Lemma \ref{sphere_chak_pal_alg} are satisfied with $ q_n^{\pm} $'s replaced by $ r_n^{\pm} $'s and hence it is easy to see that there is a $C^*$-homomorphism from $\tilde{\clq}_0$ to $\tilde{\clq}$ sending $ y, r^{+}_{n} $ to $ y_0 $ and $ q^{+}_{n} $ respectively,  which is surjective by Proposition \ref{sphere_chak_pal_description of Q'+} and is a CQG morphism by Lemma \ref{sphere_chak_pal_coproduct}. In other words, $(\tilde{\clq}_0, \Delta_0, V)$ is indeed a universal object in ${\bf Q}^\prime(D)$. It is clear that the maximal Woronowicz subalgebra of $\widetilde{\clq}_0$ for which the action is faithful, i.e. $QISO^+(D)$, is generated by $r^+_{n-1}r^{+^*}_n$, $n \geq 1$ and $r_0^+y r_1^{+^*}$, so again isomorphic with $C^{*} ( \IZ_{2}  \ast \IZ^{\infty} ).$

 This is summarized in the following:

   \bthm
   
   The universal object in the category ${\bf Q}^\prime(D)$, i.e.  $ \widetilde{QISO^{+}}( D ) $ exists and is isomorphic with $ C^{*} ( \IZ_{2}  \ast \IZ^{\infty} ).$  Moreover, the quantum group $QISO^+(D)$ is again isomorphic with $C^{*} ( \IZ_{2}  \ast \IZ^{\infty} ).$
   
   \ethm
   
 
 
   

\brmrk
This example shows that $QISO^+$ in general may not be matrix quantum group, i.e. may not have a finite dimensional fundamental unitary, even if the underlying spectral triple is of compact type. This is somewhat against the intuition derived from the classical situation, since for a classical compact Riemannian manifold, the group of isometries is always a compact Lie group, hence has an embedding into the group of orthogonal matrices of some finite dimension. 
\ermrk

\section{$QISO^+$  not having $C^*$-action}
For the rest of the article, we shall denote ${\rm ad}_V$ by $\alpha$, where $V$ is the unitary representation of $\widetilde{QISO}^+(D)$ as in the previous section.
We shall now prove that $ \alpha $ is not a  $ C^* $ action. Before that, we recall some useful properties of the so-called Toeplitz algebra from \cite{davidson}.

\bppsn

\label{sphere_toeplitz_properties}

Let $\tau_1$ be the unilateral shift operator on $l^2( \IN )$ defined by $\tau_1 ( e_n ) = e_{n - 1}, ~ n \geq 1, ~ \tau ( e_0 ) = 0.$ Then the  $C^*$ algebra $C^*(\tau_1)$ generated by $\tau_1$  on $l^2(\IN ) $, which is called the Toeplitz algebra, is a simple $C^*$ algebra and  it contains all compact operators. Moreover, the commutator of any two elements of $ C^* ( \tau_1 ) $ is compact. 

\eppsn

Let $\tau$ be the operator on $\clh$ defined by $\tau = \tau_1 \otimes {\rm id}.$ 

\blmma

\label{sphere_chak_pal_not_action_1}
One has the following polar decomposition of $B$:
$$ B = \tau \left| B \right| .$$

\elmma

{\it Proof :} We note that \bean \lefteqn{ \left| B \right| (  \left ( \begin {array} {cccc}
   e_n \\ 0 \end {array} \right ) )
= {( A - A^2 + c I )}^{\frac{1}{2}} (  \left ( \begin {array} {cccc}
   e_n \\ 0 \end {array} \right ) )}\\
 &=& \sqrt{\lambda_+ \mu^{2n} - \lambda^2_{+} \mu^{4n} + c }  (  \left ( \begin {array} {cccc}
   e_n \\ 0 \end {array} \right ) )
 = {c_+ ( n )}^{\frac{1}{2}} \left ( \begin {array} {cccc}
   e_n \\ 0 \end {array} \right ) \eean   
   and  hence $ \tau \left| B \right| (  \left ( \begin {array} {cccc}
   e_n \\ 0 \end {array} \right ) ) = {c_+ (n)}^{\frac{1}{2}} (  \left ( \begin {array} {cccc}
   e_{n - 1} \\ 0 \end {array} \right ) ) = B (  \left ( \begin {array} {cccc}
   e_n \\ 0 \end {array} \right ) ).$

Similarly, $ \tau \left| B \right| (  \left ( \begin {array} {cccc}
   0 \\ e_n \end {array} \right ) ) = B (  \left ( \begin {array} {cccc}
   0 \\ e_n \end {array} \right ) ). $ This completes the proof of the lemma. \qed

\blmma

\label{sphere_chak_pal_not_action_2}

 $$ \alpha(\tau)=\sum_{n \geq 1} \tau (P_n +Q_n) \otimes r^+_{n-1}{r^+_n}^*,$$ where $r_n^{\pm}$ are the elements of $\widetilde{QISO^+}(D)$
 introduced before. 

\elmma 

{\it Proof :} For all $ n \geq 1, $ we have  \bean \lefteqn{ \alpha ( \tau ) ( \left ( \begin {array} {cc}
   e_n \\ 0 \end {array} \right ) )  =  \widetilde{V} ( \tau \otimes {\rm id} ) {\widetilde{V}}^* \left ( \begin {array} {cccc}
   e_n \\ 0 \end {array} \right )}\\
&=& \frac{1}{2} \widetilde{V} ( \tau \otimes {\rm id} ) {\widetilde{V}}^* [ \left ( \begin {array} {cccc}
   e_n \\ e_n \end {array} \right ) + \left ( \begin {array} {cccc}
   e_n \\ - e_n \end {array} \right )  ]\\
&=& \frac{1}{2} \widetilde{V} ( \tau \otimes {\rm id} ) [ \left ( \begin {array} {cccc}
   e_n \\ e_n \end {array} \right ) \otimes  r^{+*}_n + \left ( \begin {array} {cccc}
   e_n \\ -  e_n \end {array} \right ) \otimes r^{-*}_n ]\\
&=&  \frac{1}{2} \widetilde{V} ( \tau \otimes {\rm id} ) [ \left ( \begin {array} {cccc}
   e_n \\ 0 \end {array} \right ) \otimes ( r^{+*}_n + r^{-*}_n ) + \left ( \begin {array} {cccc}
   0 \\ e_n \end {array} \right ) \otimes ( r^{+*}_n - r^{-*}_n ) ]\\
&=&  \frac{1}{2} \widetilde{V} [ \left ( \begin {array} {cccc}
   e_{n - 1} \\ 0 \end {array} \right )  \otimes ( r^{+*}_n + r^{-*}_n ) + \left ( \begin {array} {cccc}
   0 \\ e_{n - 1} \end {array} \right ) \otimes ( r^{+*}_n - r^{-*}_n )   ]\\
&=& \frac{1}{4} \widetilde{V} [ \left ( \begin {array} {cccc}
   e_{n - 1} \\ e_{n - 1} \end {array} \right ) \otimes ( 2 r^{+*}_n ) + \left ( \begin {array} {cccc}
   e_{n - 1} \\ - e_{n - 1} \end {array} \right ) \otimes 2 r^{-*}_n ]\\
 &=& \frac{1}{2} [ \left ( \begin {array} {cccc}
   e_{n - 1} \\ e_{n - 1} \end {array} \right ) \otimes r^+_{n - 1} r^{+*}_n + \left ( \begin {array} {cccc}
   e_{n - 1} \\ - e_{n - 1} \end {array} \right ) \otimes r^{-}_{n - 1} r^{-*}_n  ] = \left ( \begin {array} {cccc}
   e_{n - 1} \\ 0 \end {array} \right ) \otimes  r^+_{n - 1} r^{+*}_n. \eean
Similarly, $ \alpha ( \tau ) \left ( \begin {array} {cccc}
   0 \\  e_{n} \end {array} \right ) = \left ( \begin {array} {cccc}
   0 \\  e_{n - 1} \end {array} \right ) \otimes r^{+}_{n - 1} r^{+*}_n $ for all $ n \geq 1. $
   
   Moreover, $ \alpha ( \tau ) \left ( \begin {array} {cccc}
   e_0 \\  0 \end {array} \right ) = \alpha ( \tau ) \left ( \begin {array} {cccc}
   0 \\  e_{0} \end {array} \right ) = 0.   $  

Thus, $ \alpha ( \tau ) = \sum_{n \geq 1} \tau P_n \otimes r^+_{n - 1} r^{+*}_n + \sum_{n \geq 1} \tau Q_n \otimes r^{+}_{n - 1} r^{+*}_n = \sum_{n \geq 1} \tau ( P_n + Q_n ) \otimes r^+_{n - 1} r^{+*}_n.$
 \qed

\bthm
 \label{no_action}
 The $\ast$-homomorphism $ \alpha $ is not a $C^*$ action.

\ethm
{\it Proof :}  We begin with the observation that each of the $C^*$ algebras $\pi_{\pm}(S^2_{\mu, c})$ is nothing but the Toeplitz algebra. For example, consider $\clc:=\pi_+(S^2_{\mu, c}).$
 Clearly, $T=\pi_+(B)$ in an invertible operator with the polar decomposition given by, $ T=\tau_1 |T|,$ hence $\tau_1$ belongs to $ \clc.$ 
Thus, $\clc$ contains the Toeplitz algebra $C^*(\tau_1)$, which by Proposition \ref{sphere_toeplitz_properties} contains all compact operators.  
In particular, $C^*(\tau_1)$ must contain $\pi_+(A)$. Since $C^*(\tau_1)$ is clearly unital,  it  contains    $|\pi_+(B)|=\sqrt{\pi_+(A)-\pi_+(A)^2+cI}$, and hence also  $\pi_+(B)=\tau_1 |\pi_+(B)|$. Thus, it must contain the whole of $\clc$. 
Similar arguments will work for $\pi_-(S^2_{\mu,c}).$

Thus,  $\tau= \tau_1 \oplus \tau_1 = \pi(B)| \pi(B)|^{-1}  $ belongs to $ \pi(S^2_{\mu,c}).$ If $\alpha$ is a $C^*$ action, then for an arbitrary state $\phi$ on $QISO^+(D)$ we must have $\alpha_\phi(\tau) \equiv ({\rm id} \otimes \phi) \circ \alpha (\tau) $ is in $ \pi(S^2_{\mu,c}),$  hence $\alpha_\phi(\tau)P_+$ must belong to $\clc=\pi_+(S^2_{\mu,c})$, where $P_+$ denotes the projection 
onto $\clh_+$. By Proposition \ref{sphere_toeplitz_properties}, this implies that $[\alpha_\phi(\tau)P_+, \tau_1]$ must be a compact operator.
 We claim that for suitably chosen $\phi$, this compactness condition is violated, which will complete the proof of the theorem.

To this end, fix any nonzero real number $\theta$ and consider the sequence $\lambda_n=e^{2 \pi i n \theta}$ of complex number of unit modulus. 
 Let us consider the unital $\ast$-homomorphisms $\phi_n$, $n=0, 1, \ldots,$ where $\phi_0: C(Z_2) \raro \IC$ sends   the generator $y$ of $C(Z_2) $ to $1$ and $\phi_n$ is defined on the $n$-th copy of  $C(\IT)$ such that   $ \phi_n( r^+_{n - 1} r^{+*}_{n})=1 $.  Thus,   we have a unique unital $\ast$-homomorphism $\phi : QISO^+(D) =C(Z_2) \ast C(\IT)^{*^\infty} \rightarrow \IC$ which extends the above mentioned homomorphisms. Clearly, $ \alpha_\phi(\tau)=\sum_n \lambda_n \tau (P_n +Q_n).$ Moreover,  we see that 
\bean \lefteqn{[\alpha_\phi(\tau)P_+, \tau_1]\left ( \begin {array} {cccc}
   e_{n} \\ 0 \end {array} \right )}\\
 &=& ({\rm id} \otimes \phi ) \alpha ( \tau )\left ( \begin {array} {cccc}
   e_{n - 1} \\ 0 \end {array} \right ) - \tau ( {\rm id} \otimes \phi ) [ \left ( \begin {array} {cccc}
   e_{n - 1} \\ 0 \end {array} \right ) \otimes r^+_{n - 1} r^{+*}_n ]\\
  &=& (\lambda_{n-1}-\lambda_n) \left ( \begin {array} {cccc}
   e_{n - 2} \\ 0 \end {array} \right ) ,~~n \geq 2.\eean
Similarly, $ [\alpha_\phi(\tau)P_+, \tau_1] \left ( \begin {array} {cccc}
   0 \\ e_n \end {array} \right ) = (\lambda_{n - 1} - \lambda_n) \left ( \begin {array} {cccc}
   0 \\ e_{n - 2} \end {array} \right ) .$
Hence, the above commutator cannot be compact since $\lambda_n-\lambda_{n-1}$ does not go to $0$ as $n \rightarrow \infty.$ 
\qed
   
\bcrlre

\label{sphere_chak_pal_non_existence_wang_univ}

The subcategory  of ${\bf Q}^\prime(D)$ consisting of objects $(\tilde{\clq}, U)$ where ${\rm ad}_U$ is a $C^*$ action does not have a universal object.

\ecrlre
 {\it Proof:} By Theorem \ref{no_action}, the proof will be complete if we can show that if a universal object exists for the subcategory (say ${\bf Q}_1^\prime$) mentioned 
above, then it must be isomorphic with $\widetilde{QISO^+}(D)$. For this, consider the quantum subgroups $\widetilde{\clq}_N$, $N=1,2,...,$ 
of $\widetilde{QISO^+}(D)$ generated by $r_n^+, n=1,...,N$ and $y$. Let $\pi_N : \widetilde{QISO^+}(D) \raro \widetilde{\clq}_N$ be 
the CQG morphism given by $\pi_N(y)=y,~ \pi_N(r_n^+)=r_n^+$ for $n \leq N$ and $\pi_N(r^+_n)=1$ for $n > N.$ 

We claim that $(\widetilde{\clq}_N, U_N:=({\rm id} \ot \pi_N)\circ V)$ is an object in ${\bf Q}_1^\prime$ 
( where $V$ denotes the unitary representation of $\widetilde{QISO^+}(D)$ on $\clh$ ). 
To see this, we first note that for all $ N,$

$ ( {\rm id} \otimes \pi_N ) \alpha ( A ) = \sum^N_{n = 0} A P_n \otimes \frac{1}{2 \lambda_+} \{ \lambda_+ ( 1 + r^+_n y r^{+*}_n ) + \lambda_{-} ( 1 - r^+_n y r^{+*}_n ) \} + \sum^N_{n = 0} A Q_n \otimes \frac{1}{2 \lambda_{-}} \{ \lambda_+ ( 1 - r^+_n y r^{+*}_n ) + \lambda_{-} ( 1 + r^{+}_n y r^{+*}_n ) \} + \sum^{\infty}_{n = N + 1} A P_n \otimes \frac{1}{2 \lambda_+} \{ \lambda_+ (  1 + y ) + \lambda_{-} ( 1 - y ) \} + \sum^{\infty}_{n = N + 1} A Q_n \otimes  \frac{1}{2 \lambda_{-}} \{  \lambda_+ ( 1 - y ) + \lambda_{-} ( 1 + y ) \}.$

Among the four summands, the first two clearly belong to $ S^2_{\mu,c} \otimes \widetilde{\clq}_N.$ Moreover, the sum of the third and the fourth summand equals $ A ( 1 - \sum^{N}_{n = 1} P_n ) 
\otimes \frac{1}{2 \lambda_+} \{ \lambda_+ (  1 + y ) + \lambda_{-} ( 1 - y ) \} + A ( 1 - \sum^N_{n = 1} Q_n ) \otimes \frac{1}{2 \lambda_{-}} \{ \lambda_+ ( 1 - y ) + \lambda_{-} ( 1 + y ) \} $ which is an element of $ S^2_{\mu,c} \otimes \widetilde{\clq}_N.$

 We proceed similarly in the case of  $ B, $ to note that it is enough to show that for all $ N,$
$$  \sum^{\infty}_{n = N + 2} B P_n \otimes \frac{{c_{+}( n )}^{\frac{1}{2}} + {c_{-}( n )}^{\frac{1}{2}} }{2 c_+ ( n )}  + \sum^{\infty}_{n = N + 2} B P_n \otimes \frac{( {c_{+}( n )}^{\frac{1}{2}} - {c_{-}( n )}^{\frac{1}{2}} ) y }{2 c_+ ( n )} $$
$$ + \sum^{\infty}_{n = N + 2} B Q_n \otimes \frac{{c_{+}( n )}^{\frac{1}{2}} + {c_{-}( n )}^{\frac{1}{2}} }{2 c_- ( n )} - 
\sum^{\infty}_{n = N + 2} B Q_n \otimes \frac{( {c_{+}( n )}^{\frac{1}{2}} - {c_{-}( n )}^{\frac{1}{2}} ) y }{2 c_- ( n )} $$
 belongs to $ S^2_{\mu,c} \otimes \widetilde{\clq_{N}}.$ The norm of the second and the fourth term can be easily seen to be finite. 
The first term equals $ \frac{1}{2} B ( 1  - \sum^{N + 1}_{n = 1} P_n ) [ {( A - A^2 + c I )}^{- \frac{1}{2}} + {( A - A^2 + c I  )}^{ - 1} 
{\{ \frac{\lambda_{-}}{\lambda_+} A - {( \frac{\lambda_{-}}{\lambda_{+}} A  )}^2 + c I \}}^{\frac{1}{2}}  ] \otimes 1 $ and
 therefore belongs to $ S^2_{\mu,c} \otimes \widetilde{\clq_{N}}.$ The third term can be treated similarly.

 Thus, there is surjective  CQG morphism 
$\psi_N $ from the universal object, say $\widetilde{\clg}$, of ${\bf Q}_1^\prime$ to $\widetilde{\clq}_N$. Clearly, $ (  \widetilde{\clq}_N )_{N \geq 1}$ 
form an inductive system of objects in ${\bf Q}^\prime(D)$, with the inductive limit being $\widetilde{QISO^+}(D)$, and the morphisms 
$\psi_N$ induce a surjective CQG morphism (say $\psi$) from $\widetilde{\clg}$ to $\widetilde{QISO^+}(D)$. But $\widetilde{\clg}$ is an object
 in ${\bf Q}^\prime(D)$, so must be a sub-object (in particular, quantum subgroup) of the universal object in this category, that is, $\widetilde{QISO^+}(D)$. 
This gives a  morphism from $\widetilde{QISO^+}(D)$ onto $\widetilde{\clg}$, which is obviously the inverse of $\psi$, and hence we get the 
desired isomorphism between $\widetilde{\clg}$ and $\widetilde{QISO^+}(D).$
 \qed

\brmrk
 Since it is clear from the constructions that $QISO^+(D)$ is $\tau_I$-invariant, one has $QISO^+(D)=QISO^+_I(D)$, and hence the spectral triple considered by us also gives an example for which $QISO^+_I(D)$ does not have a $C^*$-action. Moreover, the subcategory of $Q^\prime_I$ consisting of the CQG's having $C^*$-action, does not have a universal object.
 \ermrk

Jyotishman Bhowmick : Stat-Math Unit, Indian Statistical Institute, 203, B. T. Road, Kolkata 700 108
E-mail address: jyotish r @isical.ac.in ( The support from National Board of Higher Mathematics, India,
 is gratefully acknowledged )
 
 \vspace{2mm}

Debashish Goswami ( communicating author ): Stat-Math Unit, Indian Statistical Institute, 203, B. T. Road, Kolkata 700 108
E-mail address: goswamid@isical.ac.in (  Partially supported by a project on `Noncommutative Geometry and Quantum Groups' funded by Indian National Science Academy  )

\end{document}